\documentclass[a4paper,12pt,twoside]{article}

\usepackage[austrian,american]{babel}
\usepackage[latin1]{inputenc}
\usepackage[intlimits,leqno]{amsmath}
\usepackage{amssymb}
\usepackage{amsthm}
\usepackage{graphicx,color}

\theoremstyle{plain}
\newtheorem{*theo}{*Theorem}
\newtheorem{theo}{Theorem}

\newtheorem{lemm}{Lemma}
\newtheorem{rema}{Remark}
\newtheorem{defi}{Definition}

\newcommand\I{I}
\renewcommand{\i}{{\mathrm{i}}}
\newcommand\R{{\mathbb{R}}}
\newcommand\N{{\mathbb{N}}}
\newcommand\C{{\mathbb{C}}}
\newcommand{\F}{\mathcal{F}}

\newcommand{\x}{\mathbf{x}}
\newcommand{\y}{\mathbf{y}}

\newcommand\n{{\mathbf{n}}}

\newcommand\vv{{\mathbf{v}}}

\newcommand\vO{{\mathbf{0}}}

\usepackage{makeidx}
\usepackage{epsfig}

\usepackage{latexsym}

\usepackage{pstricks,pst-node,amsfonts}

\renewcommand{\d}{{\rm d}}

\newcommand\CH{{\mathbb{H}}}
\renewcommand{\S}{\mathcal{S}}
\newcommand{\supp}{\mbox{supp}}

\begin{document}

\title{Integral equation models for thermoacoustic imaging of dissipative tissue}

\author{Richard Kowar\\ Department of Mathematics, University of Innbruck, \\
Technikerstrasse 21a/2, A-6020, Innsbruck, Austria}

\maketitle

%\abstract{
\begin{abstract}
In case of non-dissipative tissue the inverse problem of thermoacoustic imaging basically 
consists of two inverse problems. 
First, a function $\phi$ depending on the \emph{electromagnetic absorption function}, 
is estimated from one of three types of projections (spherical, circular or planar) and secondly, 
the \emph{electromagnetic absorption function} is estimated from $\phi$. 
In case of dissipative tissue, it is no longer possible to calculate explicitly the projection of $\phi$  
from the respective pressure data (measured by point, planar or line detectors). 
The goal of this paper is to derive for each of the three types of pressure data, 
an integral equation that allows estimating the respective projection of $\phi$.
The advantage of this approach is that all known reconstruction formulas for $\phi$ from the 
respective projection can be exploited.
\end{abstract}

\section{Introduction}
\label{sec:intro}

The goal of thermoacoustic imaging is to estimate the \emph{electromagnetic absorption function} of soft 
tissue so that the strong contrast between the electromagnetic absorption of cancer and normal 
tissue can be exploited (cf.~\cite{KuWanStoWan04,KucKun08,PatSch07,SchGroLenGraHal09,Wan08,XuWan06}). 
Another advantage of thermoacoustic imaging is that the technological advance allows measuring different 
types of data (cf.~\cite{Tam86,HalSchBurPal04,BurBauGruHalPal07,XuWan06}) and therefore various 
explicit reconstruction formulas can be applied (cf.~\cite{FinPatRak04,HriKucNgu08,KucKun07,KucKun08,XuFenWan02,XuXuWan02,XuXuWan03,XuWan05}). 
In order to increase the resolution of thermoacoustic imaging, we derive integral equation models that 
take dissipation into account (cf.~\cite{BurGruHalNusPal07,PatGre06,RivZhaAna06}) and allow the application 
of all known reconstruction formulas. 
In the following we explain this more precisely in case of point detectors.

First we discuss shortly the direct problem. The direct problem models the propagation of a pressure wave 
in dissipative tissue generated by very short heating of the tissue due to a laser. 
In case of homogeneous and isotropic dissipation, the direct problem can be modeled as follows 
(cf. Theorem 4.2.1 in~\cite{Hoe03})
\begin{equation}\label{modelG0}
\begin{aligned}
      p=  G*_{\x,t} f  \qquad\quad \mbox{($*_{\x,t}$ space-time-convolution)}\,,  
\end{aligned}
\end{equation}
where $G$ corresponds to an attenuated spherical wave with origin in $(\x,t)=(\vO,0)$ and $f$ 
corresponds to the  source term of $p$. In the absence of attenuation, $G$ is the retarted Green 
function of the standard wave with source term $f$.  
If $G$ is modeled appropriately, then $f$ is the same sorce term as in the absence of 
dissipation (cf. Section~\ref{sec:dp} and~\cite{SchGroLenGraHal09}), i.e. 
\begin{equation}\label{introf}
         f(\x,t)=\frac{\partial \I(t)}{\partial t}\,\phi(\x)  
\qquad\quad \mbox{($\I$ short time impulse)}\,,
\end{equation}
where $\phi$ depends on the electromagnetic absorption function, say $\phi=F(\alpha_{el})$.  
The essential requirement for $G$ is causality, i.e. the speed of the wave front (cf.~\cite{KowSchBon09})
$$
     \{(\x,T(\x))\in \R^4\,|\, T(\x):=\sup \{t>0\,|\,G(\x,\tau)=0 \mbox{ for all } \tau\leq t\} \,
$$
has to be bounded from above. Here $T(|\x|)$ denotes the travel time of the wave front.  
Below we will see that causality plays an important role in thermoacoustic imaging.

The respective inverse problem can be formulated as follows. 
Let $\Gamma$ denote a surface covered with point detectors surrounding the tissue.  
From~(\ref{modelG0}) and~(\ref{introf}), we obtain the integral equation 
\begin{equation}\label{inteqivp0}
         \int_{\R^3} \tilde G(\x-\y,t)\, \phi(\y)\,\d \y  = p(\x,t)
\qquad\mbox{ for} \quad \x\in\Gamma,\, t\in [0,T]\,,
\end{equation}
where $\tilde G :=  G *_t \frac{\partial \I(t)}{\partial t}$ and $T$ is sufficiently large. 
The inverse problem corresponds to the estimation of $\phi$ in integral equation~(\ref{inteqivp0}) 
and the subsequent inversion of $\phi=F(\alpha_{el})$.

In this paper we show that integral equation~(\ref{inteqivp0}) can be reformulated to
\begin{equation}\label{inteqivp1}
\begin{aligned}
   \int_0^\infty N(t,t')\, R_{sp}(\phi)(\x,t')\,\d t'  =  p(\x,t) 
   \qquad (N\in C^\infty (\R^2\backslash \{(0,0)\})\,,
\end{aligned}
\end{equation}
where $R_{sp}(\phi)$ denotes the \emph{spherical projection} of $\phi$, i.e. 
\begin{equation}\label{defsp}
\begin{aligned}
    R_{sp}(\phi)(\x,t) := \int_{\partial B_t(\x)} \phi(\x')\,\d \lambda^2(\x')
    \qquad\quad \x\in\R^3, t\geq 0\,.
\end{aligned}
\end{equation}
Here $\d \lambda^2$ denotes the Lebesgue measure on $\R^2$. 
The causality of $G$ ensures that $N(t,t')=0$ if $t'>t$, i.e. the upper integration limit 
in~(\ref{inteqivp1}) can be replaced by $t$. This means that the pressure function does not depend 
on information from the future.  
For some surfaces $\Gamma$, $\phi$ can be reconstructed from $R_{sp}(\phi)$ via exact or approximate 
\emph{explicit reconstruction formulas}. 
For example, if $\Gamma$ is a sphere, then~(cf. e.g.~\cite{KucKun07,KucKun08})
\begin{equation}\label{recRphip0}
\begin{aligned}
     \phi(\x) 
        = -\frac{1}{2\,\pi}\,\mbox{div}\,\int_{S^2} \n(\x')
             \left[  \frac{\partial }{\partial t} 
                        \left(   t\, R_{sp}(\phi)(\x',t)  \right)
                  \right]_{t=|\x'-\x|}
                        \d \lambda^2(\x')\,,
\end{aligned}
\end{equation}
where $\n$ denotes the exterior normal of $S^2$. 
In the absence of attenuation, if $\I(t)=\delta(t)$, then the projection $R_{sp}(\phi)$ can be 
calculated explicitly from (cf.~\cite{Joh82})
\begin{equation}\label{Rphip0}
\begin{aligned}
     \frac{\partial }{\partial t}\left( \frac{R_{sp}(\phi)(\x,t)}{4\,\pi\,t} \right) = p_0(\x,t) 
\qquad\mbox{ ($p_0$ unattenuated data)}\,.
\end{aligned}
\end{equation}
However, in the presence of attenuation, a set of $1d-$integral equations like~(\ref{inteqivp1}) 
has to be solved to obtain the projection $R_{sp}(\phi)$. 
Moreover, we see that if an efficient and accurate explicit reconstruction formula for $\phi$ from data 
$R_{sp}(\phi)$ exists, then the iterative estimation of $\phi$ via integral equation~(\ref{inteqivp0}) 
is unfavorable. 
The goal of this paper is to derive $1d-$integral equations that relate spherical, circular and planar 
projections to attenuated pressure data measured by point, line  and planar detectors, respectively. 
We note that this is a theoretical paper that is not concerned with numerical experiments.\\

The outline of this paper is as follows. 
The direct problem of thermoacoustic imaging of dissipative tissue is modeled in Section~\ref{sec:dp}.
In Section~\ref{sec:intp0} an integral equation model is derived that allows the estimation of the 
unattenuated pressure data in case $\I(t)=\delta(t)$. Then, from this model, the basic theorem of 
thermoacoustic imaging of dissipative tissue is derived in case $\I(t)\not=\delta(t)$. 
Subsequently, in Section~\ref{sec:intmod}, integral equation models for the three types of projections are 
derived that allow applying the explicit reconstruction formulas of thermoacoustic tomography. 
In the appendix it is proven that the wave equation modeled in Section~\ref{sec:dp} has a unique Green 
function satisfying the ``initial condition'' $G|_{t<0}=0$ and causality, defined as in Section~\ref{sec:dp}.

\section{Modeling of the direct problem}
\label{sec:dp}

In this section we model the direct problem of thermoacoustic imaging of dissipative tissue. 
First we summarize essential facts about wave attenuation and causality and then 
model the stress tensor and the source term for thermoacoustic imaging of dissipative tissue.
For more details we refer to~\cite{KowSchBon09,NacSmiWaa90,PatGre06,SusCob04,Sza94,Sza95}. \\

\subsection*{Causal wave attenuation in tissue}

Consider a viscous medium that is homogeneous with respect to \emph{density}, \emph{compressibility} and 
\emph{attenuation}. Then pressure waves obeying a complex attenuation law 
$\alpha_*=\alpha_*(\omega)$ satisfy the wave equation (cf.~\cite{KowSchBon09})
\begin{equation}\label{waveeq10}
\begin{aligned}
   \nabla^2 p_\gamma(\x,t)  
    -\left( D_* + \frac{1}{c_0}\,\frac{\partial}{\partial t}\right)^2 p_\gamma(\x,t)
     =  - f(\x,t)\,,
\end{aligned}
\end{equation}
where $c_0>0$ is a constant and $D_*$ is a time convolution operator with kernel defined by
\begin{equation*}
\begin{aligned}
             \hat K_*(\omega) := \frac{\alpha_*(\omega)}{\sqrt{2\,\pi}}
  \qquad \qquad \mbox{($\hat \;$ denotes Fourier transform)}\,.
\end{aligned}
\end{equation*}
The functions $\alpha_*:\R\to\C$ is called \emph{complex attenuation law} and 
$\alpha:=\mbox{Re}(\alpha_*)$ is called \emph{(real) attenuation law}. In order that $p_\gamma$ is real 
valued $\alpha$ must be even and $\mbox{Im}(\alpha_*)$ must be odd. Attenuation occurs only if $\alpha$ 
is positive. 
Moreover, the function $\alpha_*$ is restrained by the requirement of a bounded wave front speed of the 
Green function $G_\gamma$ of~(\ref{waveeq10}) (\emph{causality}). 
If the wave front speed this is bounded from above by $c_1\in (0,\infty)$, then this is equivalent to
\begin{equation*}
\begin{aligned}
    G_\gamma\left(\x,t+T(|\x|)\right) = 0\quad 
             \mbox{ if } \quad t<0 \quad \mbox{ for all } \x\not=0\,,
\end{aligned}
\end{equation*}
where $T(|\x|)\in (|\x|/c_1,\infty)$ is the travel time of the wave front from point $\vO$ to point $\x$. 
In case of a constant wave front speed $c_0$, causality is equivalent to
\begin{equation}\label{causcond}
\begin{aligned}
     K(\x,t) := 4\,\pi\,|\x|\, G_\gamma\left(\x,t+\frac{|\x|}{c_0}\right) = 0\quad 
             \mbox{ if } \quad t<0 \quad \mbox{ for all } \x\not=0\,.
\end{aligned}
\end{equation}
since $T(|\x|)=|\x|/c_0$. We note that $\hat K(\x,\omega) =  e^{-\alpha_*(\omega)\,|\x|}/\sqrt{2\,\pi}$ 
(cf. Theorem~\ref{th:unique} in the appendix).

\begin{rema}
In the literature often a less strong definition of causality is used that demands the existence 
of a (retarted) Green function that vanishes for $t<0$.  
However, this requirement is not related to the speed of the wave front. \\
\end{rema}

According to experiments the real attenuation law of a variety of viscous media similar to tissue satisfy 
(at least approximately) a frequency power law (cf.~\cite{Sza95,Web00}). 
\begin{equation*}
\begin{aligned}
    \alpha(\omega) = \alpha_0\,|\omega|^\gamma
    \qquad\mbox{ for }\qquad \gamma\in (1,2]\,,\quad \mbox{$\alpha_0=$const.}
\end{aligned}
\end{equation*}
This led to the complex attenuation laws (cf.~\cite{SusCob04,Sza94,Sza95,WatHugBraMil00})
\begin{equation}\label{powlaw1}
\begin{aligned}
    \alpha_*(\omega) = \tilde\alpha_0\,(-\i\,\omega)^\gamma
\qquad\mbox{ with }\qquad \tilde \alpha_0:= \frac{\alpha_0}{\cos(\pi\,\gamma/2)}\,,
\end{aligned}
\end{equation}
which violate causality for the range $\gamma\in (1,2]$ relevant for thermoacoustic imaging 
(cf. appendix). 
For thermoacoustic imaging we propose the following complex attenuation laws 
\begin{equation}\label{powlaw2} 
\begin{aligned}
   \alpha_*(\omega) 
       = \frac{\alpha_0\,(-\i\,\omega)}{c_0\,\sqrt{1+(-\i\,\tau_0\,\omega)^{\gamma-1}}} 
\qquad\qquad (\gamma\in (1,2],\,\tau_0>0)\,,
\end{aligned}
\end{equation}
where the square root of $1+(-\i\,\tau_0\,\omega)^{\gamma-1}$ is understood as the root that guarantees 
a nonnegative real part of $\alpha_*$.

\begin{rema}\label{rema:powfunc}
We note that the above complex power function is defined by
\begin{equation}\label{defpower}
        w^\gamma= e^{\gamma\,\left(\mbox{log}(r)+\i\,\varphi\right)}  
\qquad \mbox{ for } \qquad  
        w=r\,e^{\i\,\varphi}\in C^-\,,
\end{equation}
where $C^-:=\C\backslash\{z\in\C\,|\,\mbox{Re}(z)\leq 0,\, \mbox{Im}(z)=0\}$. 
\end{rema}

\begin{figure}[!ht]
\begin{center}
\includegraphics[height=5.0cm,angle=0]{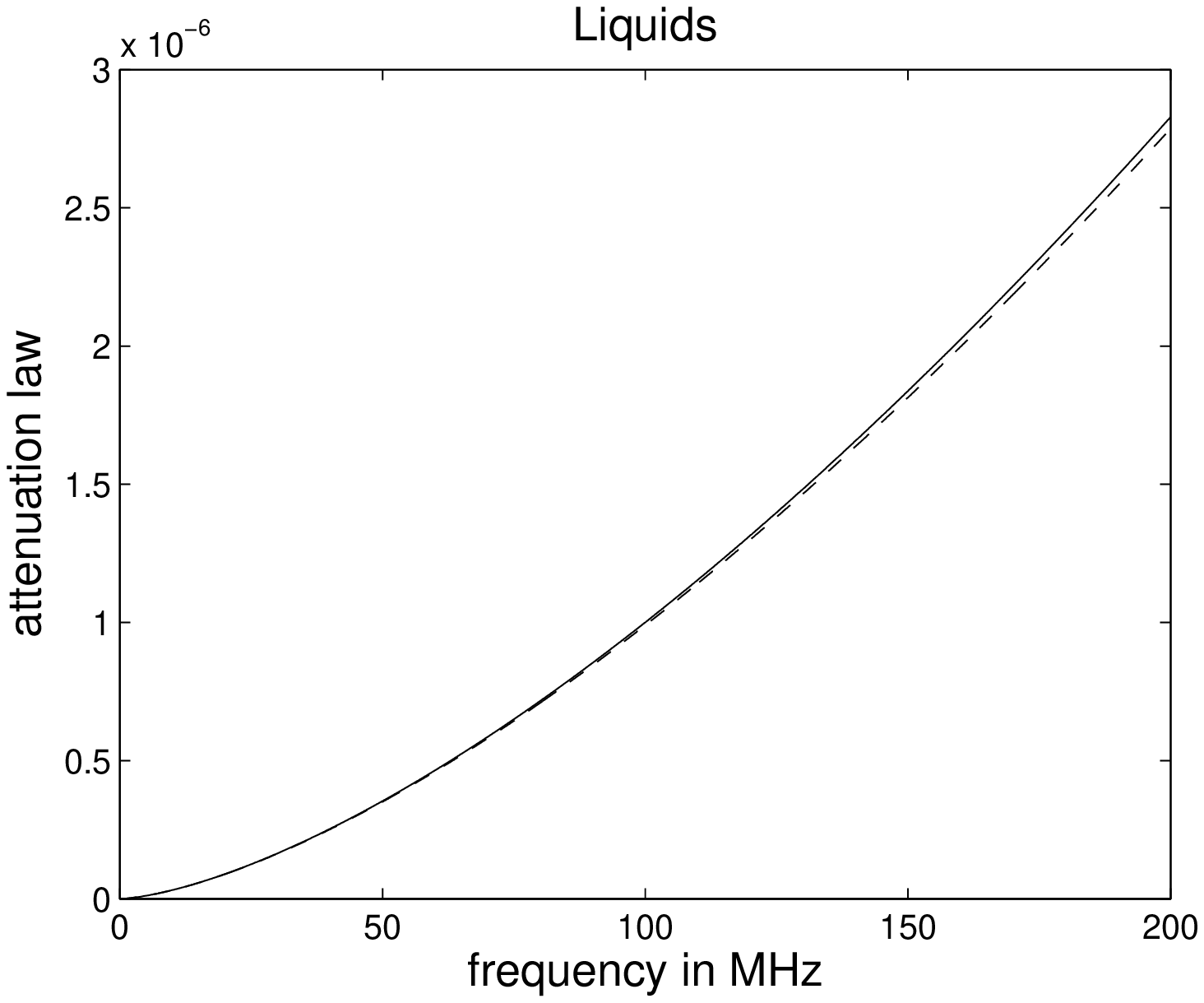}
\includegraphics[height=5.0cm,angle=0]{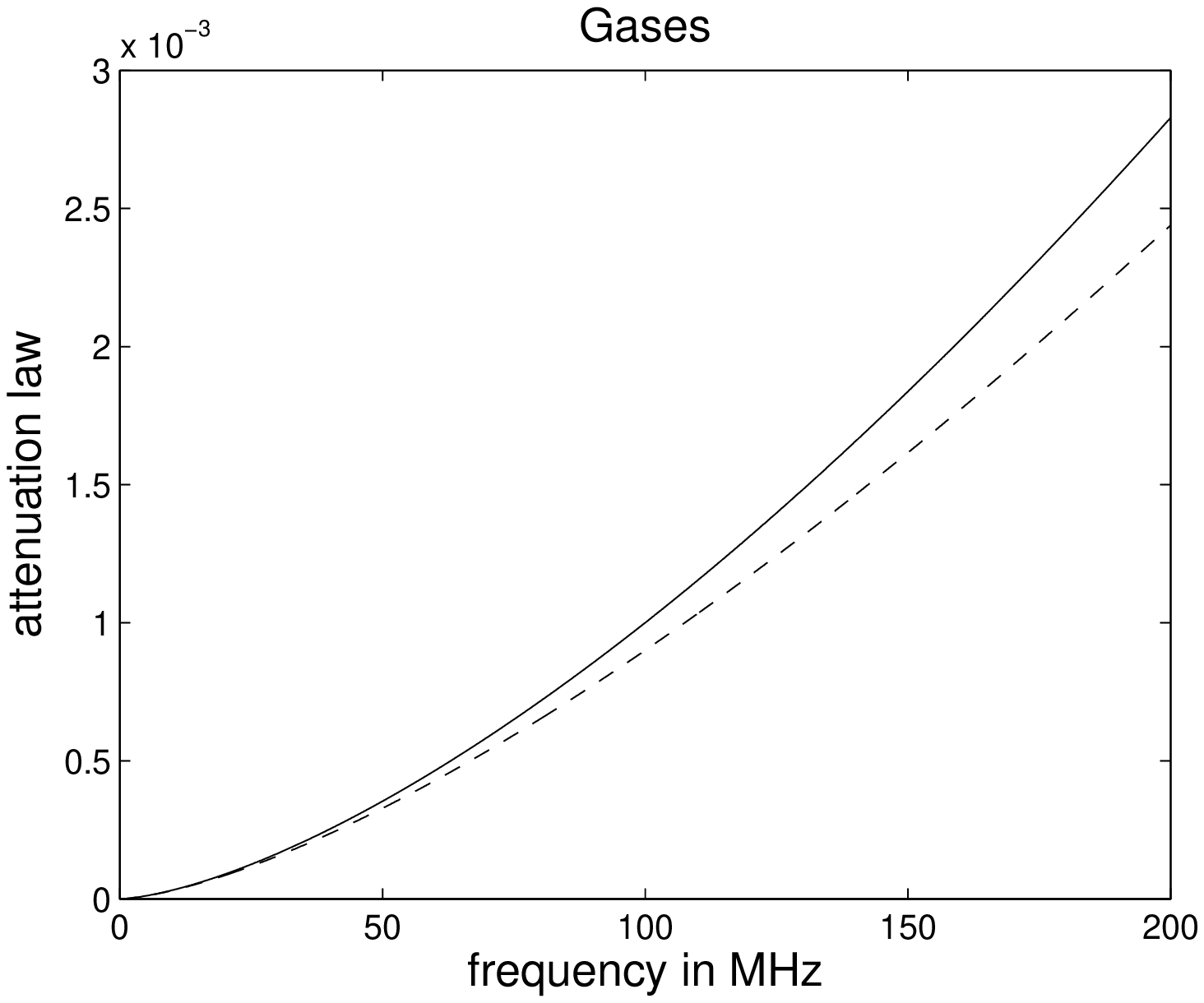}
\end{center}
\caption{Comparison of the real part of~(\ref{powlaw2}) with 
$\alpha_0:=\frac{2\,c_0\,\tau_0}{|cos(\pi\,\gamma/2)|}$ (dashed line) and the power law 
$\alpha(\omega) = |\tau_0\,\omega|^\gamma$ for $\gamma=1.5$. 
For liquids: $\tau_0=10^{-6}\,MHz$ (left picture) and for gases: $\tau_0=10^{-4}\,MHz$ (right picture) 
(cf.~\cite{KinFreCopSan00}). Experimental demonstrations of the power law are performed for the range 
$0-60\,MHz$ (cf. e.g.~\cite{Sza95}). 
}\label{fig:comp}
\end{figure}

\begin{rema}
In the appendix we prove that the Green function $G_\gamma$ of~(\ref{waveeq10}) with $\alpha_*$ 
defined as in~(\ref{powlaw1}) does not satisfy the causality requirement~(\ref{causcond}), but if 
$\alpha_*$ is defined as in~(\ref{powlaw2}) then~(\ref{causcond}) is satisfied (cf.~\cite{KowSchBon09}). 
Models~(\ref{powlaw2}) are appropriate for thermoacoustic imaging, since causality is satisfied and 
$\mbox{Re}(\alpha_*)$ is approximately a frequency power law with exponent $\gamma\in (1,2]$ for 
small frequencies (cf. Figure~\ref{fig:comp}). 
If $\gamma=2$ then $\mbox{Re}(\alpha_*)$ is the thermo-viscous attenuation law.
\end{rema}

\begin{rema}
All mathematical results in this paper are valid if $\alpha_*$ defined as in~(\ref{powlaw2}) is replaced 
by a complex attenuation law that satisfies~(\ref{causcond}) and has a monotonic increasing real part.
\end{rema}

\subsection*{Modeling of the stress tensor and the source term}

Since equation~(\ref{waveeq10}) is an integro-differential equation, it is not self-evident that 
the source term $f$ has the same structure as in the absence of wave attenuation. 
In order to model the source term we first model the temperature dependent stress tensor. 

According to~\cite{LanLif91} the stress tensor of an invisid fluid can be modeled by 
$$
    \sigma_{i,j} = - c_0^2\,\rho_0\,\alpha_{th}\,(T-T_0)\,\delta_{i,j} 
                 -p_1\,\delta_{i,j} 
\qquad\mbox{ where }\qquad p_1 = p_0 + c_0^2\,(\rho-\rho_0)\,.
$$
Here $p_0$, $\rho_0$, $c_0$, $T_0$ are reference values of the \emph{pressure}, the \emph{density}, 
the \emph{sound speed} and the \emph{temperature}, respectively, and $\alpha_{th}$ is the 
\emph{thermal expansion coefficient}. 
Wave equation~(\ref{waveeq10}) shows that a dissipative medium has a long memory (cf.~\cite{DauLio92_1}) 
and therefore we claim that attenuation is accounted for if $\sigma_{i,j}$ is replaced by 
$K_1 *_t \sigma_{i,j}$, where $*_t$ denotes time-convolution and $K_1$ is an appropriate kernel. 
We show that this ansatz implies wave equation~(\ref{waveeq10}) with the same thermoacoustic 
source term as in the absence of wave attenuation. 
(This calculation can be considered as another justification of model equation~(\ref{waveeq10}).) 
Since $\sigma_{1,1}=\sigma_{2,2}=\sigma_{3,3}$ correspond to the \emph{negative} 
total pressure in the absence of dissipation, our claim implies
\begin{equation}\label{modelp}
\begin{aligned}
    p = \tilde p_0 + c_0^2\,\,K_1 *_t(\rho-\rho_0)   + c_0^2\,\rho_0\,\alpha_{th}\,K_1 *_t (T-T_0)\,
\end{aligned}
\end{equation}
with $\tilde p_0 := p_0\, \int_\R K_1(t)\,\d t = \sqrt{2\,\pi}\,p_0\,\hat K_1(0+)$.

If the activation time of the laser is very 
short, say $\Delta t$, then the absorbed heat content of $\Omega$ at $t=\Delta t$ can be modeled 
by (cf.~\cite{GusKar93})
\begin{equation*}
\begin{aligned}
     \Delta Q(\x,\Delta t) 
        &= \alpha_{el}\,\epsilon_0(\x)\,\int e^{-\alpha_{el}\,|\x-\x_S|}\,\d \x_s \,\I(t)\,\Delta t\\
        &=:\alpha_{el}\,I_0(\alpha_{el},\x) \,\I(t)\,\Delta t,
\end{aligned}
\end{equation*}
where $\alpha_{el}$, $\epsilon_0$ and $\x_s$ denote the \emph{electromagnetic absorption coefficient}, 
the \emph{absorbed energy flow density} and the origin of an electromagnetic point source, respectively. 
Moreover, $\I=\I(t)$ denotes a short time pulse. 
On the other hand the heat content of the mass $m_\Omega$ within $\Omega$ due to the rise of temperature is 
\begin{equation*}
\begin{aligned}
     \Delta Q(\x,\Delta t) =  m_\Omega\,c_p\,(T(\x,\Delta t)-T_0),
\end{aligned}
\end{equation*}
where $c_p$ denotes the \emph{specific heat capacity} of $m_\Omega$ at constant pressure. In summary we get 
\begin{equation}\label{Q3}
\begin{aligned}
   \,\frac{\d T }{d t} 
           = \frac{\alpha_{el}\,I_0(\alpha_{el},\cdot)}{ m_\Omega\,c_p}\,\I .
\end{aligned}
\end{equation}
Since a state equation refers always to the local comoving system, the total time derivative $\d /\d t$ 
appears.
From~(\ref{modelp}) we infer the \emph{state equation}
\begin{equation}\label{stateeq}
\begin{aligned}
    \frac{1}{c_0^2}\,\frac{\d p}{\d t} 
         = K_1 *_t \frac{\d \rho}{\d t}  + \rho_0\,\alpha_{th}\,K_1 *_t \frac{\d T}{\d t}\,.
\end{aligned}
\end{equation}
% !!!
If the sound speed is not to large, then $\d /\d t$ can be replaced by $\partial / \partial t$. In the 
following we assume that this simplification is appropriate. 
% !!!
The linearized \emph{equation of motion} and the linearized \emph{continuity equation}
\begin{equation*}
\begin{aligned}
    \rho_0\,\frac{\partial \vv}{\partial t} 
          = - \nabla p 
   \qquad\mbox{ and }\qquad
\frac{\partial \rho}{\partial t} + \rho_0\,\nabla\cdot \vv = 0 
\end{aligned}
\end{equation*}
imply
$$
     \nabla^2 p - \frac{\partial^2 \rho}{\partial t^2} = 0 \,.
$$ 
From this,~(\ref{stateeq}) and~(\ref{Q3}) we infer
\begin{equation}\label{waveeq2}
\begin{aligned}
   K_1(t) *_t  \nabla^2 p(\x,t)  
    -\frac{1}{c_0^2}\,\frac{\partial^2 p(\x,t)}{\partial t^2} 
     =  - K_1(t) *_t f(\x,t) \,,
\end{aligned}
\end{equation}
with
\begin{equation*}
\begin{aligned}
      f = \frac{\alpha_{th}\,\alpha_{el}\,I_0(\alpha_{el},\cdot)}{|\Omega|\,c_p}\,
               \frac{\partial \I}{\partial t}
\qquad\quad \mbox{($|\Omega|$ volume of $\Omega$)}\,.
\end{aligned}
\end{equation*}
Comparison of~(\ref{waveeq10}) and~(\ref{waveeq2}) in Fourier space, shows that both equations 
are equivalent if
\begin{equation}\label{defK1}
\begin{aligned}
   \hat K_1(\omega) 
        := \frac{1}{\sqrt{2\,\pi}}\,\left(\frac{\omega}{c_0\,k(\omega)}\right)^2\,,
\qquad\qquad  k(\omega) := \i\,\alpha_*(\omega) + \frac{\omega}{c_0}\,.
\end{aligned}
\end{equation}
Now we can formulate the direct problem.

\subsection*{The direct problem}

Let $T,\,T_1>0$ with $T_1<<T$ and e.g. $\phi\in L^2(\Omega)$ and $\I\in \mathcal{E}(\R)$ with 
$\mbox{supp}(\I)\subset [0,T_1]$.  
The direct problem of thermoacoustic tomography is to solve wave equation~(\ref{waveeq2}) on 
$\R^3\times (0,T)$ with 
\begin{equation}\label{sourcef2}
\begin{aligned}
   f(\x,t) := \phi(\x)\,\I'(t) \qquad\mbox{ where }\qquad  
  \phi := \frac{\alpha_{th}\,\alpha_{el}}{|\Omega|\,c_p}\,I_0
\end{aligned}
\end{equation}
such that
\begin{equation}\label{initp}
     \left. p_\gamma \right|_{t<0} =0  \,.
\end{equation} 
Let $\delta(t)$ denote the \emph{delta distribution}. Frequently, it is assumed that $\I(t) =\delta(t)$, 
since the speed of light is much larger than the speed of sound. Then we have  
\begin{equation}\label{sourcef1}
\begin{aligned}
   f(\x,t) = \phi(\x)\,\delta'(t) \,
\end{aligned}
\end{equation}
and $\phi$ corresponds to an initial value function.

\begin{rema}
In the appendix we show that wave equation~(\ref{waveeq2}) has a unique Green function $G_\gamma$ 
($f(\x,t)=\delta(\x)\,\delta(t)$) satisfying $G_\gamma|_{t<0}=0$.
\end{rema}

\section{Derivation of the basic theorem}
\label{sec:intp0}

Before we explain the purpose of this section we introduce some important notions and assumptions.
\begin{itemize}
\item [(A1)] Without loss of generality we assume $c_0=1$. $\alpha_*$ is defined as in~(\ref{powlaw2}) 
             if $\gamma\in (1,2]$ and $\alpha_*:=0$ if $\gamma=0$. 
             Moreover, $k$ is defined as in~(\ref{defK1}). 
\item [(A2)] $p_\gamma$ is the unique solution of~(\ref{waveeq2}) with~(\ref{initp}) and 
             $f(\x,t) = \phi(\x)\,\I'(t)$ and 
             $\tilde p_0$ is the unique solution of~(\ref{waveeq2}) with~(\ref{initp}) 
             and $f(\x,t)= \phi(\x)\,\delta'(t)$.
\end{itemize} 
If $\I(t)\not=\delta(t)$, then we assume that the signal satisfies
\begin{itemize}
\item [(A3)]  $\I\in C(\R)$ and $\mbox{supp} (\I)\subset [0,T_1]$ for some $T_1>0$.
\end{itemize}
If $\I(t)=\delta(t)$ then $\tilde p_0=p_0$. The uniqueness of the Green function of~(\ref{waveeq2}) 
is proven in Theorem~\ref{th:unique} in the appendix.

The goal of this section is to derive an integral equation that allows the estimation of $\tilde p_0$ 
from the data $p_\gamma$. This will be the basic theorem of thermoacoustic imaging of dissipative tissue.\\

The following two lemmas provide important properties of a distribution that plays an important role in 
deriving the basic theorem.

\begin{lemm}\label{lemm01}
Let (A1) be satisfied and
\begin{equation}\label{defM}
\begin{aligned}
  \hat M_\gamma(\omega,t') 
     &:= \frac{1}{\sqrt{2\,\pi}}\,\frac{\omega}{k(\omega)}\,
             e^{\i\,k(\omega)\,|t'|} 
  \qquad (\omega,\,t'\in\R)\,.
\end{aligned}
\end{equation}
a) $M_\gamma$ satisfies for $t'\geq 0$
\begin{equation}\label{eqM}
\begin{aligned}
   K_1 *_t   \frac{\partial^2 M_\gamma}{\partial t'^2} 
    - \frac{\partial^2 M_\gamma}{\partial t^2} 
     =  0  
\quad \mbox{ with }  \quad
      \left.\frac{\partial M_\gamma}{\partial t'}\right|_{t'=0}=-\delta'(t)\,
\end{aligned}
\end{equation}
b) If $\alpha_*\not=0$ then $M_\gamma\in C^\infty(\R^2\backslash \{(0,0)\})$.\\
c) If $\alpha_*=0$ then $M_\gamma(t,t')=\delta(t-|t'|)$. 
\end{lemm}

\begin{proof}
a) In Fourier space, equation~(\ref{eqM}) and definition~(\ref{defM}) read as follows:
\begin{equation*}
\begin{aligned}
   \frac{(-\i\,\omega)^2}{(\i\,k(\omega))^2}\,
                   \frac{\partial^2 \hat M_\gamma(\omega,t')}{\partial t'^2} 
    - (-\i\,\omega)^2\,\hat M_\gamma(\omega,t') 
     =  0  \,,
\end{aligned}
\end{equation*}
\begin{equation*}
\begin{aligned}
  \frac{\partial^2 \hat M_\gamma(\omega,t')}{\partial t'^2} 
      = (\i\,k(\omega))^2\,\hat M_\gamma(\omega,t') \qquad t'>0\,.
\end{aligned}
\end{equation*}
From both identities, it follows that $M_\gamma(t,t')$ satisfies the first identity in~(\ref{eqM}). 
The second identity in~(\ref{eqM}) follows from~(\ref{defM}):
$$
    \left.\frac{\partial M_\gamma(t,t')}{\partial t'} \right|_{t'=0}
    = \F^{-1}\left\{ \frac{\i\,\omega}{\sqrt{2\,\pi}}\, \right\}(t)
    =-\delta'(t)\,.
$$
b) Let $\alpha_*\not=0$. For fixed $|t'|=T>0$, $M_\gamma$ can be considered as oscillatory integral 
(cf. Section 7.8 in~\cite{Hoe03})
\begin{equation*}
\begin{aligned}
   I_{\phi,a}(u) =  \int e^{\i\,\phi(t,\omega)} a(t,\omega)\, u(t) \d t\,\d \omega
  \qquad u\in C^\infty_0(\R)\,
\end{aligned}
\end{equation*}
with phase function $\phi(t,\omega)=-t\,\omega$ in $\Gamma=\R\times (\R\backslash \{0\})$ and 
$a(t,\omega)=\frac{1}{\sqrt{2\,\pi}}\,\frac{\omega}{k(\omega)}\,e^{\i\,k(\omega)\,T}$, since $a$ is $C^\infty$ 
and rapidly decreasing with respect to $\omega$. From Theorem 7.8.3 in~\cite{Hoe03} and
$$
   \frac{\d\phi}{\d \omega}(t,\omega)=0   \qquad\mbox{ iff }\qquad  t=0\,,
$$
it follows that $\mbox{sing\,supp\,}(I_{\phi,a}) \subseteq \{0\}$. That is to say $M_\gamma(t,T)$ is $C^\infty$ 
for $t>0$. 
Similarly, if $t=T$ is fixed then it follows that $M_\gamma(T,t')$ is $C^\infty$ for $|t'|>0$. 
This shows that $M_\gamma$ is $C^\infty$ on $\R^2\backslash \{(0,0)\}$. \\
c) The claim follows from $k(\omega)=\omega$ ($\alpha_*=0$) and property 
$\F\{ \delta(t-t')\}(\omega) = e^{\i\,\omega\,t'}/\sqrt{2\,\pi}$. 
This concludes the proof.
\end{proof}

\begin{lemm}\label{lemm02}
Let (A1) be satisfied with $\alpha_*\not=0$ and $M_\gamma$ be defined as in~(\ref{defM}). \\
a) If $t'=0$ then
\begin{equation*}
\begin{aligned}
    M_\gamma(t,0)=0 \qquad \mbox{ for } \quad t <0\,.
\end{aligned}
\end{equation*}
b) If $t'\not=0$ then $t\mapsto M_\gamma(t,t')$ is a rapidly decreasing 
$C^\infty-$function on $\R$ and satisfies
\begin{equation}\label{propM1}
\begin{aligned}
    M_\gamma(t,t')=0 \qquad \mbox{ for } \quad t\leq |t'|\,.
\end{aligned}
\end{equation}
\end{lemm}

\begin{proof}
a) From~(\ref{defM}) with $t'=0$ and~(\ref{powlaw2}) we get
\begin{equation*}
\begin{aligned}
      \sqrt{2\,\pi}\,\hat M_\gamma(\omega ,0) 
    = \frac{\omega}{k(\omega)}   
    = \frac{\sqrt{1+(-\i\,\tau_0\,\omega)^{\gamma-1}}}
           {\alpha_0+\sqrt{1+(-\i\,\tau_0\,\omega)^{\gamma-1}}}  
\qquad (\alpha_0>0)\,.
\end{aligned}
\end{equation*}
Here the power functions are defined on $\C\backslash (-\infty,0]$.  
Since $\alpha_0+\sqrt{1+(-\i\,\tau_0\,z)^{\gamma-1}}$ maps 
$$
    \{z\in\C\,|\,\mbox{Im}(z)\geq 0\}   \quad\mbox{ to }\quad
    \{w\in\C\,|\,\mbox{Re}(w)\geq \alpha_0\}\,,
$$
we see that $\hat M_\gamma(z ,0)$ is holomorphic for all $\mbox{Im}(z)\geq 0$. Hence conditions (C1) and (C2) 
of Theorem~\ref{th:lion} in the appendix are satisfied. 
Moreover, we see that there exists a polynomial $P$ such that
$$
   |\hat M_\gamma(z ,0)|   \leq    P(|z|)  \qquad \mbox{ for }\qquad \mbox{Im}(z) \geq 0\,.
$$
According to Theorem~\ref{th:lion} in the appendix, $M_\gamma(t ,0)$ vanishes for $t<0$. \\
b) Let $\alpha_*\not=0$. 
If $t'\not=0$ is fixed, then $M_\gamma(\omega,t')$ is rapidly decreasing and $C^\infty-$ 
and thus $M_\gamma(t,t')$ is a rapidly decreasing $C^\infty-$function on $\R$. Hence property~(\ref{propM1}) 
is satsified if 
\begin{equation}\label{MMM}
\begin{aligned}
    M_\gamma(t + |t'|,t')=0 \qquad \mbox{ for } \quad t < 0\,.
\end{aligned}
\end{equation}
Because of~(\ref{defM}) and the second definition in~(\ref{defK1}) we have
\begin{equation*}
\begin{aligned}
      M_\gamma(t + |t'|,t') 
          = \F^{-1}\left\{  \hat M_\gamma(\omega,t')\,e^{-\i\,\omega\,|t'|}  \right\}(t)
    =\F^{-1}\left\{ \frac{\omega}{k(\omega)}\,\frac{e^{-\alpha_*(\omega)\,|t'|}}{\sqrt{2\,\pi}}\,  \right\}(t) \, 
\end{aligned}
\end{equation*}
According to part a) of this proof and Theorem~\ref{thcaus01} in the appendix both functions 
$$
\F^{-1}\left\{ \frac{\omega}{k(\omega)} \right\}(t)  \quad\mbox{ and }\quad
\F^{-1}\left\{ \frac{e^{-\alpha_*(\omega)\,|t'|}}{\sqrt{2\,\pi}}\,  \right\}(t)
$$
vanish for $t<0$ and therefore their convolution vanishes for $t<0$, too. This proves property~(\ref{MMM}) and 
concludes the proof.
\end{proof}

\vspace{0.5cm}
\noindent
Let $\tau_1,\,\tau_2\in \R$. We recall that if $f_1$ and $f_2$ are two distributions with support 
in $[\tau_1,\infty)$ and $[\tau_2,\infty)$, respectively, then $f_1*f_2$ is well-defined 
and (cf. ~\cite{Hoe03})
\begin{equation}\label{propIM}
    \mbox{supp} (f_1 * f_2) 
           \subseteq  \mbox{supp} (f_1) + \mbox{supp} (f_2) \subseteq [\tau_1+\tau_2,\infty)\,.
\end{equation}
Now we derive the first integral relation that relates the attenuated pressure data with the 
unattenuated pressure data in case $\I(t)=\delta(t)$. 
Note that in this case $p_0\equiv \tilde p_0$ (cf. (A2)).

\begin{theo}\label{theo:pgp0}
Let (A1), (A2) with $\I(t)=\delta(t)$ be satisfied and $M_\gamma$ be defined as in~(\ref{defM}). 
Then $p_0$ solves the integral equation
\begin{equation}\label{linsysM}
\begin{aligned}
  p_\gamma(\x,t) 
     = \int_{0}^\infty M_\gamma(t,t')\,p_0(\x,t')\,\d t'\,.
\end{aligned}
\end{equation}
Here the upper limit of integration can be replaced by $t$.
\end{theo}

\begin{proof}
Let $t'\geq 0$. 
From properties~(\ref{eqM}),~(\ref{propM1}) and~(\ref{propIM}), it follows that 
\begin{equation}\label{help01}
\begin{aligned}
   K_1(t)*_t \frac{\partial^n M_\gamma(t,t')}{\partial t'^n}  = 0  
               \qquad \mbox{ if }\qquad t\leq  t' \;
           \mbox{ and } \;  n\in \{0,\,1\}\,.
\end{aligned}
\end{equation}
Moreover, $p_0$  satisfies ($\gamma=0$, $K_1(t) = \delta(t)$) 
\begin{equation}\label{help02}
\begin{aligned}
  &\nabla^2 p_0(\x,t') 
         = \frac{\partial^2 p_0(\x,t')}{\partial t'^2} 
  \qquad \mbox{ for $t'>0$} \,,\\
    & p_0(\cdot,0+) =\phi               
\qquad\mbox{ and }\qquad  
     \frac{\partial p_0}{\partial t} (\cdot,0+) =0\,.
\end{aligned}
\end{equation}
For convenience we introduce the notions 
\begin{equation}\label{opBox3}
\begin{aligned}
  \Box p :=  K_1(t)*_t (\nabla^2p)  
         - \frac{\partial^2 p}{\partial t^2}\,
\end{aligned}
\end{equation}
and
\begin{equation*}
\begin{aligned}
  q(\x,t) 
     := \int_{0}^\infty M_\gamma(t,t')\,p_0(\x,t')\,\d t'\,.
\end{aligned}
\end{equation*}
Kernel property~(\ref{propM1}) implies that the upper limit can be replaced by $t$. 
From~(\ref{eqM}) we get 
\begin{equation*}
\begin{aligned}
   \Box q(\x,t) 
    &= \int_{0}^\infty K_1(t) *_t \left[ M_\gamma(t,t')\,\nabla^2 p_0(\x,t')
                             - \frac{\partial M_\gamma(t,t')}{\partial t'} \,p_0(\x,t') 
                                    \right]\,\d t'\,,
\end{aligned}
\end{equation*}
which simplifies with integration by parts together with~(\ref{help01}) and~(\ref{help02}) to
\begin{equation}\label{Box3}
\begin{aligned}
   \Box q(\x,t) 
    &=  \phi(\x)\, \left[ K_1(t) *_t 
                     \frac{\partial M_\gamma(t,t')}{\partial t'}\right]_{t'=0} 
     = - \phi(\x)\,\,\delta(t)  *_t K_1(t) \,.
\end{aligned}
\end{equation}
Because $p_\gamma$ is the unique solution of~(\ref{waveeq2}) satisfying~(\ref{initp}) (cf. 
Theorem~\ref{th:unique} in the appendix), it follows that $p_\gamma=q$. This concludes the proof.
\end{proof}

\vspace{0.5cm}
\noindent
From Theorem~\ref{theo:pgp0}, assumption (A3) and~(\ref{propIM}),  
we get:

\begin{theo}~\label{theo:2}
Let (A1)-(A3) be satisfied and $M_\gamma$ be defined as in~(\ref{defM}). 
Then $\tilde p_0$ solves the integral equation
\begin{equation}\label{linsysM2}
\begin{aligned}
  p_\gamma(\x,t) 
     = \int_{0}^\infty  [\I(t) *_t M_\gamma(t,t')]\,\tilde p_0(\x,t')\,\d t'\,.
\end{aligned}
\end{equation}
The upper limit of integration can be replaced by $t$.
\end{theo}

\section{Integral equation models for projections of $\phi$}
\label{sec:intmod}

In this section we derive integral equation models for the three types of projections in thermoacoustic 
tomography. For various situations,  the function $\phi$ can be calculated via  
explicit reconstruction formulas from projections of $\phi$ 
(cf.~\cite{FinPatRak04,HriKucNgu08,KucKun07,KucKun08,XuFenWan02,XuXuWan02,XuXuWan03,XuWan05}). 
In case of reconstructions with limited-view, we refer to~\cite{XuWanAmbKuc03} and the literature cited 
there. \\

\noindent
In this section we use the following notions and assumptions:
\begin{itemize}
\item [1)] The region of interest (tissue) $\Omega\subset\R^3$ is a subset of the open ball 
           $B_{R_0}(\vO)$ with radius $R_0>0$.
\item [2)] $\d \lambda^m$ denotes the Lebesgue measure on $\R^m$ for $m\in\{1,\,2,\,3\}$.
\end{itemize}

\subsection*{Case of Point detectors}

The simplest setup of point detector, which guarantees stable reconstruction, is the sphere 
$\Gamma=\partial B_{R_0}(\vO)$ that encloses the region of interest. Then the set of data is
$$
          \{ p_\gamma(\x,t) \,|\, \x\in \Gamma,\,t\in [0,T] \}\,
   \quad\qquad \mbox{($T>0$ sufficiently large)\,.}
$$ 
Inserting the spherical mean representation~(\ref{Rphip0}) into integral equation~(\ref{linsysM2}) 
and performing integration by parts yield
\begin{equation*}
\begin{aligned}
  p_\gamma(\x,t) 
     = -\int_{0}^t \frac{\partial }{\partial t'}\, \left[\I(t) *_t M_\gamma(t,t')\right]\,
            \frac{R_{sp}(\phi)(\x,t')}{4\,\pi\,t'}  \,\d t'  \qquad \mbox{ for $t\geq 0$}\,,
\end{aligned}
\end{equation*}
since~(\ref{propIM}) holds and 
$$  
    \left[ \frac{R_{sp}(\phi)(\x,t')}{4\,\pi\,t'} \right]_{t'=0}
       = \left[\int_0^{t'} \tilde p_0(\x,\tau)\,\d \tau \right]_{t'=0} = 0\,.
$$
This implies

\begin{theo}\label{th:rsp}
Let (A1)-(A3) be satisfied and $M_\gamma$ be defined as in~(\ref{defM}). 
The spherical projection $R_{sp}(\phi)$ of $\phi$ satisfies 
\begin{equation}\label{linsysN1}
\begin{aligned}
  p_\gamma(\x,t) 
     = \int_{0}^t N_\gamma(t,t')\,R_{sp}(\phi)(\x,t')\, \d t'  
  \qquad \x\not=0,\, t\geq 0\,,
\end{aligned}
\end{equation}
where 
\begin{equation*}
\begin{aligned}
  N_{\gamma}(t,t')  
       &= -\frac{1}{4\,\pi\,t'}\,\I(t) *_t \frac{\partial M_\gamma(t,t')}{\partial t'}\,.
\end{aligned}
\end{equation*}
\end{theo}

\subsection*{Case of planar detectors}

In order to define a data setup of planar detectors we need 

\begin{defi}\label{defR3}
Let  $\n\in S^1$ and $E(\n,t)$ denote a plane normal to $\n$ with normal distance $t$ to the origin.
Let $\phi$ be integrable. We define the planar projecton of $\phi$ by 
\begin{equation*}
\begin{aligned}
  R_{pl}(\phi)(\n,t)
     := \int_{E(\n,t)} \phi(\x')\,\d \lambda^2(\x') \qquad \n\in S^1,\,t\geq 0\,.
\end{aligned}
\end{equation*}
\end{defi}

\vspace{0.5cm}
\noindent
Let $\Gamma=S^2$. 
The simplest setup of planar detectors, which guarantees stable reconstruction, is given by
$$
     \{ R_{pl}(p_\gamma(\cdot,t))(\n,\,R_0)\,|\, \n\in\Gamma,\,t\in [0,T] \}\,
   \quad\qquad \mbox{($T>0$ sufficiently large)\,.}
$$ 
According to~\cite{XuWan06,HalSchBurPal04} the identity
\begin{equation*}
\begin{aligned}
    R_{pl}(\tilde p_0(\cdot,s))(\n,R_0) = 2\,R_{pl}(\phi)(\n,R_0-s)  \qquad  \n\in S^1,\,s\geq 0\,
\end{aligned}
\end{equation*}
holds. This identity and integral equation~(\ref{linsysM2}) imply 
\begin{equation*}
\begin{aligned}
  R_{pl}(p_\gamma(\cdot,t))(\n,R_0)
     &= \int_{0}^t 2\,[\I(t) *_t M_\gamma(t,s)]\, R_{pl}(\phi)(\n,R_0-s)\,\d s  \\
     &= \int_{R_0-t}^{R_0} 2\,[\I(t) *_t M_\gamma(t,R_0-t')]\,  R_{pl}(\phi)\left(\n,t'\right)\,\d t'\,.
\end{aligned}
\end{equation*}
This leads to

\begin{theo}
Let (A1)-(A3) be satisfied and $M_\gamma$ be defined as in~(\ref{defM}). 
The planar projecton $R_{pl}$ of $\phi$ satisfies
\begin{equation}\label{linsysN3}
\begin{aligned}
  R_{pl}(p_\gamma(\cdot,t))(\n,R_0) 
     = \int_{R_0-t}^{R_0} N_\gamma(t,t')\,R_{pl}(\phi)(\n,t')\, \d t' 
    \qquad  \n\in S^1,\, t>0
\end{aligned}
\end{equation}
where 
\begin{equation}\label{defN3}
\begin{aligned}
  N_\gamma (t,t')  
       &=  2\,[\I(t) *_t M_\gamma(t,R_0-t')]\,.
\end{aligned}
\end{equation}
\end{theo}

\subsection*{Case of line detectors}

In this case the data setup is more complicated. Again we start with a definition.

\begin{defi}\label{defR2}
Let $\n\in S^1$ and $E(\n,0)$ be defined  as in Definition~\ref{defR3}. 
For each $\x\in E$, $l_\n(\x)$ denotes the 
line passing through $\x$ normal to $E$. For integrable $\phi$ we define the line 
integral operator by
\begin{equation*}
\begin{aligned}
     \Phi_\n
         :=  L_{\n}(\phi)(\x) 
         :=  \int_{l_\n(\x)} \phi(\x')\,\d \lambda^1(\x')
     \qquad \x\in E\,.
\end{aligned}
\end{equation*}
We define the \emph{circular projection} of $\Phi_\n$ by 
\begin{equation*}
\begin{aligned}
    R_{circ}(\Phi_\n)(\x,t) := \int_{\partial B_t(\x)\cap E} \Phi_\n(\x')\,\d \lambda^1(\x')\,
        \qquad \x\in E,\,t\geq 0\,.
\end{aligned}
\end{equation*}
\end{defi}

Actually we are concerned with three inverse problems. The first is concerned with the estimation 
of the circular projections $ R_{circ}(\Phi_\n)$ for each $\n\in S^1$, 
the second is concerned with the estimation of $\Phi_\n$ for each $\n\in S^1$ and the 
third is concerned with the estimation of $\phi$ from the set $\{\Phi_\n\,|\,\n\in S^1\}$. The latter problem 
corresponds to the inverson of the linear Radon transform and will not be discussed in this paper. 
For each $\n\in S^1$ let e.g. $\Gamma = \partial B_{R_0}(\vO)\cap E(\n)$, which guarantees a stable 
reconstruction. Consider the sets of data 
$$
     \mathcal{M}(\n) = \{ L_{\n}(p_\gamma(\cdot,t))(\x)\,|\, \x\in\Gamma,\,t\in [0,T] \}\,
   \qquad \mbox{($T>0$ sufficiently large)\,.}
$$ 
For each fixed $\n$ and data set $\mathcal{M}(\n)$ we derive an integral 
equation for $R_{circ}(\Phi_\n)$. 
In~\cite{BurBauGruHalPal07} it was shown that
\begin{equation}\label{propL}
\begin{aligned}
   L_{\n}(\tilde p_0(\cdot,s))(\x)
     = \frac{1}{2\,\pi}\,\frac{\partial }{\partial s} \int_{0}^s 
            \frac{R_{circ}(\Phi_\n)(\x,t')}{\sqrt{s^2-t'^2}} \,\d t'\,.
\end{aligned}
\end{equation}
Let $\tilde M_\gamma(t,s):=\I(t)*_t M_\gamma(t,s)$.  From~(\ref{linsysM2}) and~(\ref{propL}) and integration 
by parts we get
\begin{equation*}
\begin{aligned}
  L_{\n}(p_\gamma(\cdot,t))(\x)
     &= -\frac{1}{2\,\pi}\,\int_{0}^t \int_{0}^{s} \frac{\partial \tilde M_\gamma(t,s)}{\partial s}\, 
          \frac{R_{circ}(\Phi_\n)(\x,t')}{\sqrt{s^2-t'^2}} \,\d t' \,\d s\\
     &= -\frac{1}{2\,\pi}\,\int_{0}^t\int_{t'}^{t} \left[\frac{1}{\sqrt{s^2-t'^2}} \,
                         \frac{\partial \tilde M_\gamma(t,s)}{\partial s}\right]\,\d s\, 
          R_{circ}(\Phi_\n)(\x,t') \,\d t' \,,
\end{aligned}
\end{equation*}
since $M_\gamma(t,t)=0$ for $t>0$ and 
$$
   \frac{1}{2\,\pi}\,\left[\int_{0}^s 
            \frac{R_{circ}(\Phi_\n)(\x,t')}{\sqrt{s^2-t'^2}} \,\d t'\right]_{s=0}
      = \left[\int_0^s L_{\n}(\tilde p_0(\cdot,\tau))(\x)\,\d \tau \right]_{s=0} = 0\,.
$$
We infer

\begin{theo}
Let (A1)-(A3) be satisfied and $M_\gamma$ be defined as in~(\ref{defM}). 
For $\n\in S^1$ let $L_{\n}$ and $R_{circ}$ be defined as in Definition~\ref{defR2}. Then
\begin{equation}\label{linsysN2}
\begin{aligned}
  L_{\n}(p_\gamma(\cdot,t))(\x)
     = \int_{0}^t N_\gamma(t,t')\,R_{circ}(\Phi_n)(\x,t')\, \d t'  \qquad \x\in E,\,t>0
\end{aligned}
\end{equation}
with  kernel 
\begin{equation}\label{defN1}
\begin{aligned}
  N_\gamma(t,t')  
      &= -\frac{1}{2\,\pi}\,\int_{t'}^{t} \left[\frac{1}{\sqrt{s^2-t'^2}} \,
                         \I(t) *_t  \frac{\partial M_\gamma(t,s)}{\partial s}\right]\,\d s\,.
\end{aligned}
\end{equation}
\end{theo}

\section{Appendix}

In the appendix we prove that the frequency power law~(\ref{powlaw1}) for $\gamma\in (1,2]$ fails 
causality, defined as in Section~\ref{sec:dp}, and that the complex attenuation law proposed 
in~(\ref{powlaw2}) permits causality. 
Moreover, we prove that the attenuated wave equation has a unique Green function that vanish for $t<0$. 
For this purpose we need Theorem~4 in~\cite{DauLio92_5}, which is stated below (cf. Theorem~7.4.3 and the 
following remark in~\cite{Hoe03}).\\

Let $\CH:=\{z\in\C\,|\,\mbox{Im}(z)>0\}$ denote the upper open complex half plane and $\mathcal{S}'(\R)$ 
denotes the space of tempered distributions on $\R$ with range in $\C$.

\begin{theo}\label{th:lion}
A distribution $f\in\S'(\R)$ is causal, i.e. $\supp(f)\subseteq [0,\infty)$, if and only if 
\begin{itemize}
\item [(C1)] $\hat f:\R\to\R$ can be extended to a function $F:\CH\to\C$ that is holomorphic.
\item [(C2)] For all fixed $\eta>0$ and $\xi\in\R$, $F(\xi + \i\,\eta)$ considered as a distribution 
             with respect to the variable $\xi$ is tempered, and $F(\xi + \i\,\eta)$ converges 
             (in the sense of $\S'$) when $\eta\to 0$.
\item [(C3)] There exists a polynomial $P$ such that  
$$
               |F(z)| \leq P(|z|)  \qquad \mbox{for} \qquad \mbox{Im}(z)\geq \epsilon >0\,.
$$
\end{itemize}
If all three conditions are satisfied, then $F$ is the \emph{Fourier-Laplace transform} of $f$.
\end{theo}

\begin{rema}
The definition of the Fourier transform in this paper has a different sign as 
in~\cite{DauLio92_5} and thus $\CH$ is the upper half plane and not the lower half plane.
\end{rema}

\begin{theo} \label{coro:powlaw1}
For $\alpha_0>0$ and $\gamma\in (1,2]$ let $\alpha_*$ be defined as in~(\ref{powlaw1}). 
Then $\F^{-1}\{e^{-\alpha_*(\omega)\,|\x|}\}(t)$ is not a causal function. 
\end{theo}

\begin{proof}
Let $\x\in\R^3$ be arbitrary but fixed. \\
1) Since $e^{-\alpha_*(\omega)\,|\x|}$ is a rapidly decreasing function, its inverse Fourier transform is also a 
rapidly decreasing function and hence Theorem~\ref{th:lion} is applicable. 
The holomorphic power function $(-\i\,z)^\gamma$ (cf. Remark~\ref{rema:powfunc}) defined on $\C^-$ is the unique 
holomorphic extension of the function $\omega\in\R \mapsto (-\i\,\omega)^\gamma\in\C$ and hence 
$$
           F(z)=\exp\{-\alpha_*(z)\,|\x|\}  \qquad \mbox{ for }\quad  z\in \CH
$$
is the unique holomorphic  extension of $\omega\in\R \mapsto \exp\{-\alpha_*(\omega)\,|\x|\}\in\C$. \\
2) We show that condition (C3) cannot be satisfied for the sequence $(z_n)_{n\in\N}$ defined by
$$
         z_n:= \i\,n\in\CH   \qquad \mbox{ for }\qquad  n\in\N\,.
$$
From $-\i\,z_n=n$ and $ \cos(\frac{\pi}{2}\,\gamma)<0$ for $\gamma\in (1,2]$, it follows
$$
    |F(z_n)| 
      = \left|\exp\left\{ -\alpha_0\,\frac{(-\i\,z_n)^\gamma}{\cos(\frac{\pi}{2}\,\gamma)}  \right\}  \right|
      = \exp\left\{ \alpha_0\,\frac{n^\gamma}{|\cos(\frac{\pi}{2}\,\gamma)|}  \right\}  \,,
$$
which cannot be bounded by a polynomial $P(n)$.\\
\end{proof}

\begin{rema}
1) We note that the following modification of the power law 
\begin{equation*}
\begin{aligned}
  \alpha_*(\omega)
    = \tilde\alpha_0(\gamma)\, (-i\omega)^\gamma  + \alpha_1 \,(-\i\,\omega)
\qquad (\gamma\in (1,2])\,,
\end{aligned}
\end{equation*}
which also appears in the literature, does not satisfies causality, since
$$
    |F(z_n)| 
      = \exp\left\{ \alpha_0\,\frac{n^\gamma}{|\cos(\frac{\pi}{2}\,\gamma)|}  - \alpha_1\,n \right\} \,
$$
cannot be bounded by a polynomial for $\gamma >1$.  \\
2) Let $\gamma\in (0,1)$. It is easy to see that $|F(z)|$ is bounded by a constant, since 
$\mbox{Re}\left(\alpha_*(z)\right)$ is always positive and thus $|F(z)|$ decreases exponentially for 
$|z|\to \infty$. \\
\end{rema}

\begin{theo} \label{thcaus01}
For $\alpha_0,\,\tau_0>0$ and $\gamma\in (1,2]$, let $\alpha_*$ be defined as in~(\ref{powlaw2}). 
Then $\F^{-1}\{e^{-\alpha_*(\omega)\,|\x|}\}(t)$ is a causal function. 
\end{theo}

\begin{proof}
Let $\x\in\R^3$ be arbitrary but fixed. \\
1) Since $e^{-\alpha_*(\omega)\,|\x|}$ is a rapidly decreasing function, its inverse Fourier transform is also a 
rapidly decreasing function and thus Theorem~\ref{th:lion} is applicable. 
Let $M_{-\epsilon}:=\{z\in\C\,|\, \mbox{Im}(z)> -\epsilon\}$. 
Since $\alpha_*(z)$ is holomorphic on $M_{-\epsilon}$ for sufficiently small $\epsilon$, 
this function is the unique holomorphic extension of $\omega\in\R \mapsto\alpha_*(\omega)\in\C$. 
Thus 
$$
           F(z)=\exp\{-\alpha_*(z)\,|\x|\}  \qquad \mbox{ for }\quad  z\in \CH
$$
is the unique extension of $\omega\in\R \mapsto \exp\{-\alpha_*(\omega)\,|\x|\}\in\C$. 
This proves conditions (C1) and (C2) of Theorem~\ref{th:lion}. \\
3) Since the inequality in (C3) is equivalent to
$$
     \exp\{-\mbox{Re}(\alpha_*(z))\} \leq P(|z|)  \qquad \mbox{for} \qquad \mbox{Im}(z)\geq \epsilon >0\,,
$$
(C3) is satisfied if
$$
           \mbox{Re}(\alpha_*(\CH)) \subseteq [0,\infty)\,.
$$
Since 
$$
    \mbox{Re}(\alpha_*) = \frac{\alpha_0}{c_0}\,(w_1\,z_2 + w_2\,z_1)  \qquad\quad (z_2>0)
$$
with
$$
   f(z) 
       := \frac{1}{\sqrt{1+(-\i\,\tau_0\,z)^{\gamma-1}}} =: w_1 +\i\,w_2 
      \qquad\mbox{ for }\qquad z\in\CH\,,
$$
we have to show that
\begin{equation}\label{propw}
       w_1>0 \qquad\mbox{ and }\qquad     w_2\,z_1 \geq 0\,.
\end{equation}
Let $f_1=-\i\,\tau_0\,z$, $f_2(z) := \sqrt{1+z^{\gamma-1}}$ and $f_3(z) := \frac{1}{z}$, so that 
$f=f_3\circ f_2 \circ f_1$. Moreover, let 
\begin{equation*}
\begin{aligned}
   M_1 &:= \{z_1+\i\,z_2\in\C\,|\, z_1<0\,, z_2>0 \}\,,  \\
   M_2 &:= \{z_1+\i\,z_2\in\C\,|\, z_1>0\,, z_2>0 \}\,,\\
   M_3 &:= \{z_1+\i\,z_2\in\C\,|\, z_1>0\,, z_2<0 \}\,, \\
   N_1 &:= \{\i\,z_2\in\C\,|\, z_2>0 \}\,\qquad\mbox{and}\qquad
   N_2 &:= \{z_1\in\C\,|\, z_1>0 \}\,.
\end{aligned}
\end{equation*}
We see that 
\begin{equation*}
\begin{aligned}
  &  f_1(M_1)\subseteq M_2\,,\qquad
     f_1(N_1)\subseteq N_2\,,\qquad
     f_1(M_2)\subseteq M_3\,, \\
  & f_2(M_2)\subseteq M_2\,,\qquad
    f_2(N_2)\subseteq N_2\,,\qquad
    f_2(M_3)\subseteq M_3\,.
\end{aligned}
\end{equation*}
From
$$
          \frac{1}{z_1 +\i\,z_2} = \frac{z_1 -\i\,z_2}{z_1^2 + z_2^2}  \qquad\quad z_1^2 + z_2^2\not=0\,,
$$
we get
\begin{equation*}
\begin{aligned}
    f_3(M_2)\subseteq M_3\,,\qquad
    f_3(N_2)\subseteq N_2\,,\qquad
    f_3(M_3)\subseteq M_2\, .
\end{aligned}
\end{equation*}
Therefore we end up with
\begin{equation*}
\begin{aligned}
    f(M_1)\subseteq M_3\,,\qquad
    f(N_1)\subseteq N_2\,,\qquad
    f(M_2)\subseteq M_2\, ,
\end{aligned}
\end{equation*}
since $f=f_3\circ f_2 \circ f_1$. In other words condition~(\ref{propw}) is satisfied and hence condition 
(C3) holds. This concludes the proof.
\end{proof}

\begin{theo}\label{th:unique}
For $\alpha_0,\,\tau_0>0$ and $\gamma\in (1,2]$, let $\alpha_*$ be defined as in~(\ref{powlaw2}). 
Then equation~(\ref{waveeq2}) with $f(\x,t)=\delta(\x)\,\delta(t)$ and
\begin{equation}\label{init}
    p_\gamma|_{t<0} = 0
\end{equation}
has a \emph{unique solution} 
$$
      p_\gamma(\x,t) 
              = \F^{-1}\left\{  
                 \frac{\exp\{\i\,k(\omega)\,|\x|\}}{4\,\pi\,|\x|}
                       \right\}
\quad\mbox{ with }\quad  
k(\omega) := \i\,\alpha_*(\omega) + \frac{\omega}{c_0}\,.
$$
\end{theo}

\begin{proof} 
Since $p_\gamma$ is the Green function, we write $G_\gamma$ instead of $p_\gamma$. 
The Helmholtz equation of~(\ref{waveeq2}) reads as follows
$$
            \nabla^2 \hat G_\gamma +k^2\,\hat G_\gamma = -\frac{\delta(\x)}{\sqrt{2\,\pi}}
$$
and has the only two solutions
$$
  \hat G_\gamma(\x,\omega) = \frac{1}{\sqrt{2\,\pi}}\,\frac{\exp\{s\,\i\,k(\omega)\,|\x|\}}{4\,\pi\,|\x|} 
  \qquad\quad s=\pm 1\,.
$$
If condition~(\ref{init}) holds, then according to Theorem~\ref{th:lion}, we have for $z=z_1+\i\,z_2\in\CH$ and 
$\x\not=\mathbf{0}$: 
\begin{equation}\label{conds}
          e^{-s\,\mbox{Re}(\alpha_*(z))\,|\x|}\,e^{-s\,z_2\,\frac{|\x|}{c_0}} 
            = |e^{s\,\i\,k(z)\,|\x|}|  \leq P(|z|)    \qquad  z_1\in\R,\,\,z_2>0
\end{equation}
for some polynomial $P$. From Theorem~\ref{thcaus01} and Theorem~\ref{th:lion}, we get
\begin{equation*}
          e^{-\mbox{Re}(\alpha_*(z))\,|\x|}
             \leq Q(|z|)    \qquad  z_1\in\R,\,\,z_2\geq \epsilon >0
\end{equation*}
for some polynomial $Q$, i.e.~(\ref{conds}) holds with $P:=Q$ for $s=1$. 
If $s=-1$, then the left hand side of~(\ref{conds}) grows exponentially, since 
$\mbox{Re}(\alpha_*(\omega))>0$ is increasing, i.e. condition~(\ref{init}) cannot be satisfied for $s=-1$. 
Hence there exists a unique solution. 
\end{proof}

\end{document}